\documentclass[a4paper,11pt]{amsart} 
\usepackage{amscd}             
\usepackage{xypic}  
\usepackage{amssymb}
\usepackage{amsthm}
\usepackage{epsf}
\newtheorem{thm}{Theorem}[section]   

\newtheorem{lemma}[thm]{Lemma}
\newtheorem{prop}[thm]{Proposition}

\newtheorem{rem}[thm]{Remark}

\def\im{\operatorname{Im}}

\def\c1{\operatorname{c_1}}
\def\c2{\operatorname{c_2}}
\def\Cliff{\operatorname{Cliff}}
\def\gon{\operatorname{gon}}

\def\Hilb{\operatorname{Hilb}}

\def\ZZ{{\mathbf Z}}

\def\QQ{{\mathbf Q}}
\def\PP{{\mathbf P}}

\def\P{{\mathcal P}}

\def\O{{\mathcal O}}
\def\I{{\mathcal I}}

\def\U{{\mathcal U}}
\def\V{{\mathcal V}}

\def\x{\times}                   
\def\iso{\simeq}

\def\+{\oplus}                   
\def\*{\otimes}                  
\def\hpil{\longrightarrow}       
\def\khpil{\rightarrow}

\def\Pic{\operatorname{Pic}}

\begin{document}

  \title{On two conjectures for curves on $K3$ surfaces}
  \author{Andreas Leopold Knutsen}  
  
\address{\hskip -.43cm Andreas Leopold Knutsen, Dipartimento di Matematica, Universit\`a di Roma Tre, 
Largo San Leonardo Murialdo 1, 00146, Roma, Italy. e-mail {\tt knutsen@mat.uniroma3.it}}


\thanks{Research supported by a Marie Curie Intra-European Fellowship within the 6th 
European Community Framework Programme}

\begin{abstract}
We prove that the gonality among the smooth curves in a complete linear system on a $K3$ surface is constant except for the 
Donagi-Morrison example. This was proved by Ciliberto and Pareschi \cite{cp} under the additional condition that the linear system is ample. 

As a consequence we prove that exceptional curves on $K3$ surfaces satisfy
the Eisenbud-Lange-Martens-Schreyer conjecture \cite{elms}  and explicitly describe such curves. They turn out to be natural extensions of the Eisenbud-Lange-Martens-Schreyer 
examples of exceptional curves on $K3$ surfaces.
\end{abstract}

\maketitle

\section{Introduction} \label{intro}

In connection with their work \cite{hm}, Harris and Mumford
conjectured  that the gonality should be constant among the smooth curves in a linear system on a $K3$ surface. (The conjecture is unpublished.)
Subsequently, Donagi and Morrison \cite{dm}  pointed out the following counterexample:

\vspace{0,3cm}

\noindent {\bf The Donagi-Morrison example } (cf. \cite[(2.2)]{dm}). 
 Let $\pi:S \khpil \PP^2$ be a $K3$ surface of genus $2$, i.e. a double cover of  $\PP^2$ branched along a smooth sextic, and let $L:= \pi ^* \O_{\PP^2}(3)$. The arithmetic genus of the curves in $|L|$ is $10$. The smooth curves in the codimension one linear subspace $|\pi ^* H^0\O_{\PP^2}(3)| \subset |L|$ 
are biellliptic, whence with gonality $4$. On the other hand the general 
curve in $|L|$ is isomorphic to a smooth plane sextic and therefore has 
gonality $5$. 

\vspace{0,3cm}

Ciliberto and Pareschi \cite[Thm. A]{cp} proved that this is indeed the only 
counterexample when $L$ is ample. The first aim of this note is to show that this result
holds without the ampleness assumption. 
That is, we will prove:

\vspace{0,3cm}

\noindent {\bf Theorem 1.} {\it Let $S$ be a $K3$ surface and $L$ a globally generated line bundle on $S$. If the gonality of the smooth curves in $|L|$ is not constant, then $S$ and $L$ are as in the Donagi-Morrison example.} 

\vspace{0,3cm}

It has also been known that this result would follow from the Eisenbud-Lange-Martens-Schreyer 
conjecture on {\it exceptional} curves posed in \cite[p.~175]{elms} (see \S \ref{S:p2}). (Recall that 
any smooth curve $C$ satisfies $\Cliff C +2 \leq \gon C \leq \Cliff C +3$ and the
curves for which $\gon C = \Cliff C +3$ are conjectured to be very rare and called {\it exceptional}.) 
In \cite[Thm. 4.3]{elms} an infinite series of examples of exceptional curves lying on $K3$ surfaces is constructed. 
The line bundles in these cases are not ample (cf. also 
\cite[Remark (c), p.~36]{cp}), showing that there are interesting cases appearing  when the line bundles are not ample.

We will consider a generalization of these examples:

\vspace{0,3cm}

\noindent {\bf ``Generalized ELMS examples''.} Let $L$ be a  line bundle on a $K3$ surface $S$ such that $L \sim 2D + \Gamma$ with $D$ and $\Gamma$ smooth curves  
satisfying $D^2 \geq 2$, $\Gamma^2=-2$ and $\Gamma.D=1$. Assume furthermore that there is no line bundle $B$ on $S$ satisfying $0 \leq B^2 \leq D^2-1$ and $0< B.L-B^2 \leq D^2$. 

Then $|L|$ is base point free and all the smooth curves in $|L|$ are exceptional, of genus $g = 2D^2+2 \geq 6$, Clifford index $c=D^2-1=\frac{g-4}{2}$ and Clifford dimension $r=\frac{1}{2}D^2+1$. Moreover, for any 
smooth curve $C \in |L|$ the Clifford index is computed only by $\O_C(D)$.
(Recall that the {\it Clifford dimension} of a smooth curve is the minimal value of $\dim |A|$, 
where $A$ computes the Clifford index.)
\vspace{0,3cm}

We will prove the assertions in the example in Proposition \ref{prop:elmsexa}. The examples in
\cite[Thm. 4.3]{elms} have $\Pic S \iso \ZZ[D] \+ \ZZ[\Gamma]$ with $D$ and $\Gamma$ as above, in which case
the nonexistence of a divisor $B$ satisfying the conditions above can easily be verified. 

As in \cite{elms}, the curves in the ``generalized ELMS examples'' satisfy the
Eisenbud-Lange-Martens-Schreyer conjecture.

The second main result of this note is:

\vspace{0,3cm}

\noindent {\bf Theorem 2.} {\it Let $C$ be a smooth exceptional curve on a $K3$ surface $S$.
Then $C$ is either a smooth plane sextic belonging to the Donagi-Morrison example or
$\O_S(C)$ is as in the generalized ELMS examples.

In particular, $C$ satisfies the  Eisenbud-Lange-Martens-Schreyer conjecture.}

\vspace{0,3cm}

We remark that the proof of Theorem 1, as well as the assertions in the ``generalized ELMS examples''
(in Proposition \ref{prop:elmsexa})  {\it do not} use the theorem of 
Green and Lazarsfeld \cite{gl} about constancy of the Clifford index (as in the case of Ciliberto and Pareschi's paper, cf. \cite[Rem. p.~32]{cp}). The latter enters the picture only in the proof of Theorem 2.

We prove Theorem 1 by adding a suitable deformation-degeneration argument 
to the arguments of \cite[\S 1 and \S 2]{cp}. (We do not make use 
of \cite[\S 3]{cp}.) We therefore use the same notation and conventions as in \cite{cp} and refer the reader to that paper for background material.

The note is organised as follows.

In Section \ref{S:1} we obtain sharper versions of Lemma 2.2 and Proposition 2.3 in \cite{cp} and introduce an incidence variety, slightly different from the one considered in \cite[\S 3]{cp}, that we will need in the proof of Theorem 1.

In Section \ref{S:p1} we prove Theorem 1. The idea is as follows: Since Theorem 1 holds when $L$ is ample, by \cite{cp}, the ideal way to prove it would be to deform $(S,L)$ so as to
\begin{itemize}
\item[(i)] keep the nonconstancy of the gonality among the smooth curves in $|L|$, and
\item[(ii)] make $L$ ample.
\end{itemize}

The condition (i) is easily preserved in a codimension two subspace of the moduli space: one just needs to keep the two line bundles $M$ and $N$ such that $L \sim M+N$ coming from the instability of the well-known vector bundle considered in \cite{cp}. 

Condition (ii) is not possible to achieve, but we will show that we can make $L$ ``almost ample'', in the sense that there is a unique rational curve $\Gamma$ such that $\Gamma.L=0$. Moreover, we will show that
$H:=L-\Gamma$ is globally generated and we will prove Theorem 1 by degenerating to the special curves
$C'' \cup \Gamma$ in the linear system $|L|$, with $C'' \in |H|$ smooth, and using the incidence variety from Section \ref{S:1}.

In Section \ref{S:p2} we prove the assertions in the ``Generalized ELMS examples''
in Proposition \ref{prop:elmsexa} and then we prove Theorem 2, which at this point is just a combination of Theorem 1 with the well-known theorem of Green and Lazarsfeld \cite{gl}.

\section{Some useful results} \label{S:1}

We first obtain some strengthenings of \cite[Lemma 2.2 and Prop. 2.3]{cp} in Lemma \ref{lemma:2.2} and Proposition \ref{prop:2.3}, respectively, as we will need these stronger versions in the proof of Theorem 1.

\begin{lemma} \label{lemma:2.2}
Let $L$ be a base point free line bundle on a $K3$ surface $S$ and assume that
$L \sim M+N$ with $h^0(M) \geq 2$, $h^0(N) \geq 2$, $M.N=k$ and 
$L^2 \geq 4k-4$.

Then either
\begin{itemize}
\item[(a)] there is a smooth curve in $|L|$ of gonality $\leq k$; or
\item[(b)] $M \sim N +\Gamma$ (possibly after interchanging $M$ and $N$), 
for a smooth rational curve $\Gamma$ such that $\Gamma.N=1$. In particular, 
$L^2=4k-2$.
\end{itemize}
\end{lemma}

\begin{proof}
Among all the decompositions satisfying the conditions in the lemma, we pick one for which $k$ is minimal, say $L \sim M_0+N_0$ with $M_0.N_0=k_0 \leq k$. If $k_0=k$, we let $M_0=M$ and $N_0=N$. (Note that we have $k_0 \geq 2$ as $L$ is globally generated, cf. \cite{S-D}.)

If $M_0 \sim N_0$, then $M_0$ is nef, as $L$ is. If it were not base point free, then $M_0 \sim lE+\Gamma$, for $l \geq 2$, a smooth elliptic curve $E$ and a smooth rational curve $\Gamma$ such 
that $E.\Gamma=1$, by \cite{S-D}. One then easily sees that $|E|$ induces a pencil of degree $\leq k_0$ on all the curves in $|L|$ and we are in case (a).

By symmetry we can therefore assume that $M_0.L \geq N_0.L$ and $h^0(N_0-M_0)=0$. 
We now show that either we are in case (a) or we can find a new decomposition
$L \sim M' + N'$ satisfying the following properties:
\begin{eqnarray}
\label{eq:newdec1} M' \geq M_0, \; N' \leq N_0, \; M'.N' = k_0; \\
\label{eq:newdec2} {M'}^2 \geq {N'}^2 >0; \\
\label{eq:newdec3} N' \; \mbox{is globally generated with} \; h^0(N') \geq 2;\\
\label{eq:newdec4}  h^1(M')=h^1(N')=0; \\
\label{eq:newdec5}  \mbox{the base divisor $\Delta'$ of $|M'|$ satisfies $\Delta'.L=0$.}
\end{eqnarray}

If $N$ is not nef, then there is a smooth rational curve $\Gamma$ such 
that
$\Gamma.N_0 <0$. Therefore $\Gamma.M_0 >0$ as $L$ is nef, and $h^0(M_0+\Gamma) \geq h^0(M_0) \geq 2$, $h^0(N_0-\Gamma)=h^0(N_0) \geq 2$ and  
\[ (M_0+\Gamma).(N_0-\Gamma) =k_0+\Gamma.N_0 -\Gamma.M_0 +2 \leq k_0. \]
Hence, the minimality of $k_0$ implies $\Gamma.N_0=-1$ and $\Gamma.M_0=1$, so that
$(M_0+\Gamma).(N_0-\Gamma) =k_0$. 
In particular, continuing the process, we 
reach a decomposition $L \sim M' + N'$ satisfying 
\eqref{eq:newdec1} with $N'$ nef. As above, if $N'$ is not base point free, 
then $N' \sim lE+\Gamma$, for $l \geq 2$, a smooth elliptic curve $E$ and a smooth rational curve $\Gamma$ such 
that $E.\Gamma=1$. One then easily sees that $|E|$ induces a pencil of degree 
$\leq k_0$ on all the curves in $|L|$ and we are in case (a). 
Otherwise \eqref{eq:newdec3} is satisfied.

If ${N'}^2=0$, then $N_0.L =k_0$, so that all the curves in $|L|$ would carry a pencil of degree $k_0$, and we are in case (a) again. Otherwise ${N'}^2 >0$, and as
$M'.L \geq M_0.L \geq N_0.L \geq N'.L$, we have ${M'}^2 \geq {N'}^2 >0$, so that
\eqref{eq:newdec2} is satisfied. In particular, $h^1(N')=0$. Moreover, the above argument with
$M_0$ and $N_0$ substituted by $N'$ and $M'$ respectively, shows that any $\Delta>0$ satisfying $\Delta^2=-2$ and $\Delta.M'<0$, must satisfy $\Delta.M'=-1$. Hence $h^1(M')=0$ by \cite[Thm.~1]{klvan} and \eqref{eq:newdec4} is satisfied.

Let now $\Delta'$ be the (possibly zero) base divisor of of $|M'|$
and assume that  $\Delta'.L >0$. 

If $h^1(M'-\Delta')>0$, then by \cite{S-D} we have $M' -\Delta' \sim lE$ for a smooth elliptic curve $E$ and an integer $l \geq 2$. But then $|E|$ is easily seen to induce a pencil of degree $\leq k_0$ on the 
curves in $|L|$, so that we are in case (a).

If $h^1(M'-\Delta')=0$, then $M'.\Delta'=\frac{1}{2}{\Delta'}^2<0$ by Riemann-Roch, as 
$h^0(M'-\Delta')=h^0(M')$ and $h^1(M')=0$.
Moreover, $N'.\Delta' \geq -M'.\Delta'+1$, by assumption. Hence
\[
 (M'-\Delta').(N'+\Delta') =M'.N'-\Delta'.N' + \Delta'.M' -{\Delta'}^2 < k_0, 
\]
a contradiction on the minimality of $k_0$.

Therefore, \eqref{eq:newdec5} is proved.

Now we set $R':=M'-N'$.
Then the condition  $L^2 \geq 4k_0-4$ is equivalent to 
${R'}^2 \geq -4$. We have showed above that $h^2(R')=0$.

Let now $D \in |N'|_s$ (the locus of smooth curves in $|N'|$, with notation as in \cite{cp}) 
 and consider $\O_D(M')$. 

We now claim that
\begin{eqnarray} 
\label{eq:2.2.4} \mbox{$\O_D(M')$ is base point free if and only if $R'$ is not a smooth rational curve} \\
\nonumber \mbox{satisfying $R'.N'=1$ (in which case $L^2=4k_0-2$, so that $k=k_0$).\hspace{0,8cm}}
\end{eqnarray}

If ${R'}^2 = -4$, then $L^2=4k_0-4={N'}^2+{M'}^2+2k_0$, whence 
${N'}^2+{M'}^2=2k_0-4$, and it follows that ${N'}^2 \leq k_0-2$, 
since ${N'}^2 \leq {M'}^2$. Therefore $\deg \O_D(M')=k_0 \geq {N'}^2+2 =2g(D)$ 
and $\O_D(M')$ is base point free.
  
If ${R'}^2 \geq -2$, then $R >0$ by Riemann-Roch and the fact that 
$h^2(R')=0$.

We have
  \[ \deg \O_D(M') =N'.M' = (N'+R').M'={N'}^2+R'.N'=
2g(D)-2+R'.N', \]
  so  $\O_D(M')$ is base point free if $R'.N' \geq 2$. If 
$R'.N' \leq 1$, we must have ${R'}^2 = -2$ by \cite{S-D}, as $N'$ is globally generated.
We will now show that $R'$ is irreducible with $R'.N'=1$.

We have $R'.L =2R'.N'-2$, whence $R'.N'=1$ and $R'.L=0$ by the nefness of $L$. So there has to exist a smooth rational curve $\Gamma \leq R'$ 
such that $\Gamma.N'=1$.  Now $2N'+\Gamma \leq L$, and since
$h^0(2N'+\Gamma) \geq \frac{1}{2}(2N'+\Gamma)^2+2 \geq h^0(L)$, 
we must have $R'=\Gamma$. Hence \eqref{eq:2.2.4} is proved.

By \cite[Lemma 2.2]{cp} and the conditions \eqref{eq:newdec1}-\eqref{eq:2.2.4},
we are therefore in case (a) unless $R' \sim M'-N'$ is a smooth rational curve 
and $R'.N'=1$. In this case $k_0=k$ so that $M_0=M$ and $N_0=N$. We have $(M-N)^2=-2$, so that $M -N >0$ by Riemann-Roch, and
since $M-N \leq M'-N' = R'$, we have $M=M'$ and $N=N'$ and we are in case (b).
\end{proof}

\begin{prop} \label{prop:2.3}
Keep the same hypotheses and notation as in \cite[Prop. 2.3]{cp}.

If we are in case (b) of \cite[Prop. 2.3]{cp}, then all the smooth curves in $|L|$ have gonality 
$d$ and Clifford index $d-3$, so are exceptional.

If we are in case (c) of \cite[Prop. 2.3]{cp} with $\rho(g,d,1) <0$, then the following additional conditions hold:
\begin{itemize}
\item[$(c_6)$] $M.L \geq N.L$ and $h^0(N-M)=0$ unless $M \sim N$;
\item[$(c_7)$] $M$ is not of the form $M \sim N+\Delta$, with $\Delta$ a smooth rational curve
such that $\Delta.N=1$ (and $\Delta$ is the base divisor of $|M|$);
\item[$(c_8)$] for any smooth, irreducible $D \in |N|$, we have that $\O_D(M)$ is base point free.
\end{itemize}
If, furthermore, the gonality among the smooth curves in $|L|$ is not constant, then
\begin{itemize}
\item[$(c_9)$] the general $C' \in |L|$ satisfies $\Cliff C=\Cliff C'=\Cliff \O_{C'}(N)=d-2$ and 
$\gon C'=d+1$ (whence is exceptional);
\item[$(c_{10})$] $L^2 \geq 4d-2$ and $M-N >0$;
\item[$(c_{11})$] $M^2 >0$ and $N^2 >0$.
\end{itemize}
\end{prop}

\begin{proof}
Assume we are in (b) of \cite[Prop. 2.3]{cp}. Then, for any smooth $C' \in |L|$, one easily sees that
 $\O_{C'}(N)$ contributes to the Clifford index of $C'$, as $h^0(N) =h^0(L-N) \geq 2$, so that
\begin{eqnarray*}
\Cliff C' & \leq & \Cliff \O_{C'}(N)= \deg \O_{C'}(N)-2(h^0(\O_{C'}(N))-1) \\
& \leq & L.N -2(h^0(N)-1) \leq L.N- 2(\frac{1}{2}N^2+1) \\
& = & L.N-N^2-2 =  N.(N+\Delta)-2 =c_2(E_{C,A})-1-2 \\
& = & d-3 =\gon C-3 \leq \gon C'-3.
\end{eqnarray*}
Since $\Cliff C'=\gon C'-2$ or $\gon C'-3$ by \cite[Thm. 2.3]{cm}, we must have
$\Cliff C'= \gon C'-3=\gon C-3=d-3$, 
so that all $C' \in |L|$ have the same gonality $d$ and the same Clifford index $d-3$.
Hence they are all exceptional.

Assume now that we are  in (c) of \cite[Prop. 2.3]{cp}. Note that 
$d=M.N$ and that $L^2 \geq 4d-4$ as $\rho(g,d,1) <0$.

By \cite[Lemma 2.1]{cp}, either $h^0(M-N) >0$ or the sequence in \cite[$(c_5)$ in Prop. 2.3]{cp},
\begin{equation} \label{eq:split}
0 \hpil M \hpil E_{C,A} \hpil N \hpil 0, 
\end{equation} 
splits. Hence we can without loss of 
generality assume $(c_6)$ by symmetry.

To prove $(c_7)$, assume by contradiction that $M \sim N+\Delta$, with $\Delta$ a smooth rational curve
such that $\Delta.N=1$. Then  $h^1(\Delta)=h^1(M-N)=0$ by Riemann-Roch. Hence \eqref{eq:split} splits, contradicting the fact that $E_{C,A}$ is globally generated off a finite set, as $\Delta$ is the base divisor of $|M|$
(cf. \cite[Lemma 1.1(d)]{cp}).

Next note that $(c_8)$ follows from $(c_7)$ exactly as in the proof of \eqref{eq:2.2.4} above.

Now assume that the gonality among the smooth curves in $|L|$ is not constant. Then $(c_{11})$ follows as otherwise $M$ (or $N$) would cut out on every $C' \in |L|$ a pencil of degree $ \leq d$.

As one easily sees that $\O_{C'}(N)$ contributes to the Clifford index of any $C' \in |L|$, we get
$\Cliff C' \leq \Cliff \O_{C'}(N)=d-2$, whence by \cite[Thm. 2.3]{cm},
$\gon C' \leq \Cliff C'+3 \leq d-2+3=d+1$, so that
$(c_9)$ follows. 

By \cite[Cor. 1.3 and Prop. 2.1]{elms} we have $g(C') \geq 2\Cliff C'+4 =2d$, whence
$L^2 \geq 4d-2$ and the rest of $(c_{10})$ follows using $(c_6)$ and Riemann-Roch.
\end{proof}

As the last preparatory material for the proof of Theorem 1, we will now consider an incidence variety that is slightly different from the one in \cite[\S 3]{cp}.

Assume that we are in case (c) of \cite[Prop. 2.3]{cp} with $\rho(g,d,1) <0$ (without the assumption
that the gonality is not constant). Consider the incidence $\I_{L,N,d} \subset |L| \x |N|_s \x \Hilb^{d}(S)$
defined by
\[  \I_{L,N,d} := \Big\{ (C,D,Z) \; | \; Z \subset C \; \mbox{and} \; Z \in |\O_D(M)| \Big\}, \]
and let $p^1_{L,N,d}$, $p^2_{L,N,d}$ and $p^3_{L,N,d}$ be the projections.

\begin{lemma} \label{lemma:inci}
Assume that $M \not \sim N$. Then
\begin{itemize}
\item[(a)] $\I_{L,N,d}$ is irreducible of dimension $\dim |L|+1$;
\item[(b)] the projection $\I_{L,N,d} \khpil |L| \x \Hilb^{d}(S)$ is an isomorphism onto its image;
\item[(c)] if $C \in |L|_s$ lies in $\im p^1_{L,N,d}$, then $\gon C=d$.
\end{itemize}
\end{lemma}

\begin{proof}
The Hodge index theorem, $(c_6)$ and the fact that $M \not \sim N$ imply $D^2 =N^2 < M.N=d$. Therefore, two 
distinct $D_1,D_2 \in |N|_s$ cannot share the same $Z$ and (b) follows.

Consider the incidence $\I_{N,d} \subset |N|_s \x \Hilb^{d}(S)$ given by 
$\I_{N,d}:=\{ (D,Z) \; | \; Z \in |\O_D(M)| \}$.
This is smooth, 
irreducible of dimension $\dim |N|+ \dim |\O_D(M)|= d$, using the fact that $h^1(\O_D(M))=0$ for 
reasons of degree. For any $D \in |N|_s$ and any $Z \in |\O_D(M)|$, we have 
\begin{equation} \label{eq:diml0}
\dim |L \* \I_Z|= \dim |L|-d+1 >0,
\end{equation}
as can be computed from
\begin{equation} \label{eq:suc}
0 \hpil M \hpil L \* \I_Z \hpil \omega_D \hpil 0.
\end{equation}
and the fact that $h^1(M)=0$ by property $(c_3)$ in \cite[Prop. 2.3]{cp}. Therefore
$\I_{N,d}=\im (p^2_{L,N,d} \x  p^3_{L,N,d})$ and the dimension of any fiber 
$(p^2_{L,N,d} \x  p^3_{L,N,d})^{-1}([D,Z])$ is
$\dim |L \* \I_Z|= \dim |L|-d+1$.
This proves (a) and the fact that $Z$ does not impose independent conditions on $L$
implies also (c).
\end{proof}

\section{Proof of Theorem 1} \label{S:p1}

Let $L$ be a globally generated line bundle on a $K3$ surface $S$ and assume that the gonality of the smooth curves in $|L|$ is not constant. Let $d$ be the minimal gonality among the smooth curves in $|L|$ and let $C \in |L|$ be a smooth $d$-gonal curve. Then 
$\rho(g,d,1) <0$ by Brill-Noether theory, where $g=\frac{1}{2}L^2+1$ is the genus of $C$. Hence we are in case (c) of \cite[Prop. 2.3]{cp} and 
the conditions $(c_1)$-$(c_5)$ therein and $(c_6)$-$(c_{11})$  in Proposition \ref{prop:2.3} are satisfied. In particular, we have:
\begin{eqnarray}
\label{eq:c3} L \sim M + N, \; M^2 >0, \; N^2>0, \; h^i(M)=h^i(N)=0 \; \mbox{for} \; i=1,2; \\
\label{eq:c2} \mbox{$N$ is globally generated; \hspace{3,5cm}} \\
\label{eq:c9} \mbox{the general $C' \in |L|$ satisfies 
$\Cliff C=\Cliff C'=\Cliff \O_{C'}(N)=d-2$} \\
\nonumber \mbox{and $\gon C'=d+1$ (whence is exceptional); \hspace{2cm}} \\
\label{eq:c10} L^2 \geq 4d-2 \; \mbox{and} \; M-N >0. \hspace{3cm} 
\end{eqnarray}

Assume now, to get a contradiction, 
that we are not in the Donagi-Morrison example. We claim that
\begin{eqnarray}
\label{eq:ass1} h^1(M-N) >0, \hspace{2cm} \\
\label{eq:ass2} \mbox{$M$ and $N$ are linearly independent in $\Pic S$.} 
\end{eqnarray}

Indeed, if $h^1(M-N)=0$, then \eqref{eq:split} splits, so that $E_{C,A} \iso M \+ N$ and
$h^1(E_{C,A} \* E_{C,A}^*)=0$ and we are in the Donagi-Morrison example by 
\cite[Cor. 1.6]{cp}, a contradiction. Moreover, if 
 $M$ and $N$ are linearly dependent in $\Pic S$, then $M \sim mB$ and $N \sim nB$ for a nef $B$ in $\Pic S$ and positive integers $m$ and $n$, whence the contradiction $h^1(M-N)=0$.

Let $f:\mathfrak{S} \khpil \U$ denote the Kuranishi deformation of $S=S_0$, $0 \in \U$. 
Then $\U$ is smooth of dimension $20$, cf. \cite{Kod} or \cite[VIII, Thm. 7.3]{bpv}.
Let now $\V' \subset \U$ be the submanifold to which both line bundles $L$ and $N$ lift. By \eqref{eq:c3} and \eqref{eq:ass2}, $\V'$ is smooth of dimension $18$ by \cite[Thm. 14]{Kod}. Again by \cite[Thm. 14]{Kod},
there is a Zariski-open dense subset $\V \subset \V'$ such that for any $t \in \V -\{0\}$, we have
that $S_t$ is a smooth $K3$ surface and $\Pic S_t$ has rank two, where $S_t$ denotes the surface corresponding to $t \in \V$. Therefore, letting $L_t$, $N_t$ and $M_t:=L_t-N_t$ denote the deformations of 
$L=L_0$, $N=N_0$ and $M=M_0$, we have
\begin{equation} \label{eq:picst}
\Pic _{\QQ} S_t \iso \QQ[N_t] \+ \QQ[L_t].
\end{equation}

The next lemma shows that the ``nonconstancy of gonality'' is preserved by the deformation.

\begin{lemma} \label{lemma:goncliffgen}
Let $t \in \V -\{0\}$ be general.
Then
\begin{itemize}
\item[(i)]  there is a smooth curve $C_t \in |L_t|$ with $\gon C_t=d$; 
\item[(ii)] the general $C_t' \in |L_t|$ satisfies $\Cliff C_t=\Cliff C'_t=\Cliff \O_{C'_t}(N_t)=d-2$ 
and $\gon C'_t=d+1$ (whence is exceptional).
\end{itemize}
\end{lemma}

\begin{proof}
If (i) does not hold, then, as $M_t.N_t=d$, we must have $(M_t-N_t)^2=-2$, $(M_t-N_t).L_t=0$ and 
$h^0(M_t-N_t)=1$ by Lemma \ref{lemma:2.2}. Hence also $(M-N)^2=-2$ and $(M-N).L=0$, whence
$h^0(M-N)=1$, so that $h^1(M-N)=0$ by Riemann-Roch, contradicting \eqref{eq:ass1}.

As in the proof of Proposition \ref{prop:2.3}, one sees that
$\Cliff C'_t \leq \Cliff \O_{C}(N')= d-2$ for general $C '_t \in |L_t|$,
and equality must hold and $\gon C'_t=d+1$ by \cite[Thm. 2.3]{cm} as these hold for $t=0$ by \eqref{eq:c9}, proving (ii).
\end{proof}

We will need the following technical lemma about divisors on $S_t$:

\begin{lemma} \label{lemma:nonampio}
Let $t \in \V -\{0\}$ be general. Then there is a unique smooth, rational curve $\Gamma_t \subset S_t$ such that $\Gamma_t.L_t=0$. Furthermore, 
\begin{itemize}
\item[(i)]   $\Gamma_t \sim _{\QQ} -2N_t + \frac{2(d+N^2)}{L^2}L_t$;
\item[(ii)] $M_t-\Gamma_t$ is globally generated and $(M_t-\Gamma_t)^2 >0$;
\item[(iii)] $L_t-\Gamma_t \sim N_t + (M_t-\Gamma_t)$ is the only decomposition satisfying 
$h^0(N_t) \geq 2$, $h^0(M_t-\Gamma_t) \geq 2$ and $N_t.(M_t-\Gamma_t) \leq d-1$ (in fact,
$N_t.(M_t-\Gamma_t) = d-1$);
\item[(iv)] $M_t-N_t-\Gamma_t >0$ and $h^1(M_t-N_t-\Gamma_t)=0$.
\end{itemize} 
\end{lemma}

\begin{proof}
By \cite[Thm. A]{cp} and Lemma \ref{lemma:goncliffgen} we have that $L_t$ cannot be ample, so that there is a 
smooth, rational curve $\Gamma_t \subset S_t$ such that $\Gamma_t.L_t=0$. 

As $h^1(M_t)=0$ by \eqref{eq:c3} and $N_t$ is globally generated by \eqref{eq:c2}, by \cite[Thm.~1]{klvan} we can only have
$(\Gamma_t.N_t,\Gamma_t.M_t)=(0,0)$ or $(1,-1)$. Writing $\Gamma_t \sim aN_t + bL_t$ with
$a,b \in \QQ$ by \eqref{eq:picst} we obtain $-2=\Gamma_t^2=aN_t.\Gamma_t$, whence $a=-2$,
$N_t.\Gamma_t=1$ and (i) easily follows. This also proves that $\Gamma_t$ is unique.

Note that \eqref{eq:ass2}, \eqref{eq:c10} and the Hodge index theorem imply $N^2 <d$, so that 
$\frac{2(d+N^2)}{L^2} \leq 1$ by \eqref{eq:c10}. Hence $M_t-\Gamma_t  \sim N_t + (1- \frac{2(d+N^2)}{L^2})L_t$ 
is nef. Moreover, any smooth elliptic curve $E_t \subset S_t$ satisfies 
$E_t.(M_t-\Gamma_t) \geq E_t.N_t \geq 2$ by \cite{S-D} as $N_t$ and $L_t$ are globally generated, whence (ii) follows by \cite{S-D}.

Now assume $L_t-\Gamma_t \sim A_t + B_t$ satisfies 
$h^0(A_t) \geq 2$, $h^0(B_t) \geq 2$ and $A_t.B_t \leq d-1$. We have $\Gamma_t.(A_t + B_t)=2$.
Since $A_t.(B_t+\Gamma_t) \geq d$ and $B_t.(A_t+\Gamma_t) \geq d$ by Lemma \ref{lemma:2.2}  and 
\eqref{eq:c10}, we can only have $A_t.B_t = d-1$ and
$\Gamma_t.A_t=\Gamma_t.B_t=1$. Writing $A_t \sim xN_t + yL_t$ with
$x,y \in \QQ$ by \eqref{eq:picst} we therefore obtain $x=1$. Moreover, from
\[ d=A_t.(B_t+\Gamma_t)=(N_t+yL_t).(-N_t+(1-y)L_t) \]
we obtain $2y(d+N^2)=y(1-y)L^2$. Hence either $y=0$ or $1-y=\frac{2(d+N^2)}{L^2}$ and (iii) is proved.

Finally, note that $M_t-N_t-\Gamma_t \sim _{\QQ} (1-\frac{2(d+N^2)}{L^2})L_t$. Hence
$h^1(M_t-N_t-\Gamma_t)=0$ as $L_t$ is big and nef and $\frac{2(d+N^2)}{L^2} \leq 1$. Moreover,
$h^0(M_t-N_t-\Gamma_t) >0$ by Riemann-Roch, and by Lemma \ref{lemma:goncliffgen} combined with $(c_7)$ in 
Proposition \ref{prop:2.3}, we must have $M_t-N_t-\Gamma_t >0$, proving (iv).
\end{proof}

If now $(S_t,L_t)$ is as in the Donagi-Morrison example, then $L_t \sim 3B_t$ with $B_t^2=2$, and as this 
is preserved for $t=0$, also $(S,L)$ is as in the Donagi-Morrison example, a contradiction. 

To reach the desired contradiction, thus proving Theorem 1, we can therefore assume that the following additional conditions are satisfied:
\begin{eqnarray}
\label{eq:add1} L \sim H +\Gamma, \; \mbox{with $H$ globally generated and $\Gamma$ a smooth,} \\
\nonumber    \mbox{rational curve such that $\Gamma.M=-1$ and $\Gamma.N=1$;\hspace{0,5cm}} \\ 
\label{eq:add2} \mbox{$M-\Gamma$ is globally generated and $h^1(M-\Gamma)=0$;\hspace{0,5cm}} \\
\label{eq:add3} \mbox{$H \sim N + (M-\Gamma)$ is the only decomposition satisfying 
$h^0(N) \geq 2$, \hspace{-0,5cm}} \\
\nonumber    \mbox{$h^0(M-\Gamma) \geq 2$ and $N.(M-\Gamma) \leq d-1$ (in fact,
$N.(M-\Gamma) = d-1$);\hspace{-0,5cm}} \\
\label{eq:add4} M-N-\Gamma >0 \; \mbox{and} \; h^1(M-N-\Gamma)=0.\hspace{2cm}
\end{eqnarray}

Consider now the incidence $\I_{L,N,d} \subset |L| \x |N|_s \x \Hilb^{d}(S)$
defined in Section \ref{S:1}. By Lemma \ref{lemma:inci}, we see that we would reach the desired contradiction, that is, that $\gon C'=d$ for general $C' \in |L|$, if we show that
\begin{equation} \label{eq:dimfibragen}
\dim (p^1_{L,N,d})^{-1}(C') =1 \; \mbox{for general} \; C' \in \im p^1_{L,N,d}.
\end{equation}
We show \eqref{eq:dimfibragen} by showing that 
\begin{equation} \label{eq:im}
C'' \cup \Gamma \in \im p^1_{L,N,d} \; \mbox{for general} \; C'' \in |H|
\end{equation}
and
\begin{equation} \label{eq:dimfibraspec}
\dim (p^1_{L,N,d})^{-1}(C'' \cup \Gamma) =1 \; \mbox{for general} \; C'' \in |H|.
\end{equation}
(Recall that $C''$ is smooth by \eqref{eq:add1}.)

To this end we will need:

\begin{lemma} \label{lemma:d-1}
All the smooth curves in $|H|$ have gonality $d-1$, and for the general smooth 
$C'' \in |H|$ we have
\begin{itemize}
\item[(i)]  $\dim W_{d-1}^1(C'')=0$ and 
\item[(ii)] $C''$ contains some $W \in |\O_{D}(M-\Gamma)|$ for some $D \in |N|_s$ and $|\O_{C''}(W)|$ is a $g^1_{d-1}$.
\end{itemize}
\end{lemma}

\begin{proof}
 By \cite[Prop. 2.3]{cp}, Lemma \ref{lemma:2.2} and \eqref{eq:add3} the minimal gonality of a 
smooth curve in $|H|$ is $d-1$, as $H^2 =L^2-2 \geq 4d-4=4(d-1)$ by \eqref{eq:c10}. Hence, by 
\cite[Lemma 1.2 and Cor. 1.6]{cp} the first assertion follows from \eqref{eq:add4} by using 
the vector bundle $N \+ (M-\Gamma)$.

For a general $C'' \in |H|$, let $|A''|$ be a $g^1_{d-1}$ on $C''$. Then from \eqref{eq:add3}, \eqref{eq:add4} 
and \cite[Prop. 2.3]{cp} we have $E_{C'',A''} \iso \O_S(N) \+ \O_S(M-\Gamma)$
and from property $(c_8)$ in Proposition \ref{prop:2.3} we have that $|\O_D(M-\Gamma)|$ is base point free for any
$D \in |N|_s$. Pick a $W \in |\O_D(M-\Gamma)|$. From
\[ 0 \hpil M-\Gamma \hpil H \* \I_W \hpil \omega_D \hpil 0 \]
and \eqref{eq:add2} we see that $|H \* \I_W|$ is globally generated off $W$. For general $D$ and $W$, the general member of $|H \* \I_W|$ is smooth by Bertini
(and the base point freeness of $|\O_D(M-\Gamma)|$). Moreover, one easily computes that $\dim |H \* \I_W|= \dim |H|-d+2$, so that $|\O_{C''}(W)|$ is a $g^1_{d-1}$.
Using the standard exact sequence involving $E_{C'',A''}$, 
\[ 0 \hpil H^0(A'')^* \* \O_S \hpil  \O_S(N) \+ \O_S(M-\Gamma) \hpil 
\omega_{C''}-A'' \hpil 0 \]
(cf. \cite[(2), p.~17]{cp}), one easily sees that, for any 
$W \in |A''|$, one has 
$h^0(N \* \I_W)=h^0(\O_{C''}(N)(-A''))=h^0(\O_S) \+ h^0(N-(M-\Gamma)) =1$ 
and $h^0((M-\Gamma) \* \I_W)=h^0(\O_{C''}(M-\Gamma)(-A''))
=h^0(M-\Gamma-N)+1 \geq 2$, where we have used \eqref{eq:add2} and \eqref{eq:add4}. 
Therefore, there is a $D \in |N|$ containing $W$. From what we saw above, for general $C''$ and $W$, this $D$ is smooth. Moreover, there is an $M' \in |M-\Gamma|$ containing $W$ but not $D$, so that
$W=D \cap M'$, whence $W \in |\O_D(M-\Gamma)|$. This proves (ii).

Consider the incidence $\I_{H,N,d-1}$. We have $M-\Gamma \not \sim N$ by \eqref{eq:add4} and we have just seen that $p^1_{H,N,d-1}$ is dominant, whence by Lemma \ref{lemma:inci}(i) its fibers are one-dimensional, proving (i).
\end{proof}

Now \eqref{eq:im} follows from Lemma \ref{lemma:d-1}(ii). 

Pick a general $C'' \in |H|$ satisfying \eqref{eq:im}. Then by Lemma \ref{lemma:inci}(b), we have
that \linebreak $\dim (p^1_{L,N,d})^{-1}(C'' \cup \Gamma) = \P_1 \cup \P_2$, where
\[ \P_1= \Big \{ Z \; | \; Z \in |\O_D(M)| \; \mbox{for some} \; D \in |N|_s \; \mbox{and} \; Z \subset C'' \Big\} \]
and
\[ \P_2= \Big \{ Z \; | \; Z \in |\O_D(M)| \; \mbox{for some} \; D \in |N|_s, \; Z =W \cap \{x\}, \; 
x=\Gamma \cap D  \; \mbox{and} \; W \subset C'' \Big\}. \]

As $W \in |\O_D(M-\Gamma)|$, we have $\dim \P_2=1$ by Lemma \ref{lemma:d-1}.

To compute $\dim \P_1$, consider the incidence 
 $\I \subset |H| \x |N|_s \x \Hilb^{d}(S)$
defined by
\[  \I= \Big\{ (C'',D,Z) \; | \; Z \subset C'' \; \mbox{and} \; Z \in |\O_D(M)| \Big\}, \]
and let $q_1$, $q_2$ and $q_3$  be the projections. As in Lemma \ref{lemma:inci} the projection
$\I \khpil |H| \x \Hilb^{d}(S)$ is an isomorphism onto its image, and as we can assume that $q_1$ is dominant,
we have
\begin{eqnarray*}
\dim \P_1 & = & \dim q_1^{-1}(C'')= \dim \I- \dim |H| = \dim \im (q_2 \x q_3) + \dim |H \* \I_Z| - \dim |H| \\
          & = & \dim \I_{N,d} + \dim |H \* \I_Z| - \dim |H| = d -(d-1) =1.
\end{eqnarray*}
Here $\I_{N,d}:=\{ (D,Z) \; | \; Z \in |\O_D(M)| \} \subset |N|_s \x \Hilb^{d}(S)$ 
is the incidence variety in the proof of Lemma \ref{lemma:inci} (where we showed that
$\dim \I_{N,d} =d$) and $\dim |H \* \I_Z| =  \dim |H| -(d-1)$ is easily calculated from
\eqref{eq:suc} tensored by $\O_S(-\Gamma)$, using Riemann-Roch and \eqref{eq:add2}.

Hence \eqref{eq:dimfibraspec} follows and Theorem 1 is proved.

Note that by \cite[Thm. 3.1]{cp}, we have the following consequence of Theorem 1:

\begin{thm} \label{thm:zero}
Let $S$ be a $K3$ surface and $L$ a globally generated line bundle on $S$, not as in the
Donagi-Morrison example. Let $g$ be the genus and $d$ the gonality of the smooth curves in $|L|$. 

If $\rho(d,g,1) <0$, then, for the general smooth $C \in |L|$, we have
$\dim W^1_d=0$.
\end{thm}

\section{Proof of Theorem 2} \label{S:p2}

We will first prove the assertions in the ``generalized ELMS examples''.

\begin{prop} \label{prop:elmsexa}
Let $L$ be a line bundle on a $K3$ surface $S$ such that $L \sim 2D + \Gamma$ with $D$ and $\Gamma$ smooth curves  
satisfying $D^2 \geq 2$, $\Gamma^2=-2$ and $\Gamma.D=1$. Assume furthermore that there is no line bundle $B$ on $S$ satisfying $0 \leq B^2 \leq D^2-1$ and $0< B.L-B^2 \leq D^2$. 

Then $|L|$ is base point free and all the smooth curves in $|L|$ are exceptional, of genus $g = 2D^2+2 \geq 6$, Clifford index $c=D^2-1=\frac{g-4}{2}$ and Clifford dimension $r=\frac{1}{2}D^2+1$. Moreover, for any 
smooth curve $C \in |L|$ the Clifford index is computed only by $\O_C(D)$.
\end{prop}

\begin{proof}
Since $\Gamma.L=0$ and $D.L >0$ we have that $L$ is nef. Moreover, any smooth elliptic curve $E$ on $S$ must satisfy $E.L =2E.D+E.\Gamma \geq 2E.D \geq 4$, as $D$ is nonrational, whence $|L|$ is base point free by \cite{S-D}.

Now set $k:=D^2+1=D.(D+\Gamma)$. 
For any smooth $C \in |L|$ one computes
\[ \Cliff C \leq \Cliff \O_C(D) =k-2. \]
Assume that $d:=\gon C \leq k$. Then $\rho(g,d,1) <0$, whence by \cite[Prop. 2.3]{cp} there is a globally generated $N \in \Pic S$ such that $h^0(N) \geq 2$, $h^0(L-N) \geq 2$, $h^1(N)=h^1(L-N)=0$, $N.(L-N) \leq k$ and $(L-N)^2 \geq N^2 \geq 0$ (the latter by 
$(c_6)$ in Proposition \ref{prop:2.3} and by Riemann-Roch on $N$). 

We want to show that $N \sim D$.

The Hodge index theorem and the fact that $L \not \sim 2N$ yield $N^2 \leq k-1$. If equality holds, then for the same reason we have $N.(L-N)=k$, whence $N.L=2k-1=N.(2D+\Gamma) \geq 2D.N$. It follows 
that $D.N \leq k-1$ and $N \sim D$ by the Hodge index theorem, as desired.

If $N^2 \leq k-2=D^2-1$, then the assumption on the nonexistence of $B$ implies $N.(L-N)=k$. Let now $F:=D-N$. Then one easily computes $k=D.(L-D)=F.(F+\Gamma)+N.(L-N)= F.(F+\Gamma)+k$, whence $F.(F+\Gamma)=0$. 
As $h^1(L-N)=0$  we must have $\Gamma.(L-N) \geq -1$ by \cite[Thm.~1]{klvan}, whence $\Gamma.N=0$ or $1$. As $1=\Gamma.D=\Gamma.(F+N)$, we conclude that $\Gamma.F=F^2=0$ and $\Gamma.N=1$. We then get
\begin{eqnarray*}
 F.L & = & F.(2D+\Gamma) =2D.F =2(N+F).F=2N.F =(L-2N-\Gamma).N \\
& = & N.(L-N)-N^2-\Gamma.N=k-N^2-1.
\end{eqnarray*}
But then $0< F.L \leq k-1$, a contradiction on the nonexistence of $B$.

It follows that $L \sim D + (D+\Gamma)$ is the only decomposition
satisfying $h^0(D) \geq 2$, $h^0(D+\Gamma) \geq 2$ and $D.(D+\Gamma) \leq k:=D^2+1$.
Therefore, we cannot be in case (c) of \cite[Prop. 2.3]{cp}, by condition $(c_7)$ in Proposition
\ref{prop:2.3}. Hence we must be in case (b) and by Proposition \ref{prop:2.3},  
all the smooth curves in $|L|$ have gonality 
$k+1$ and Clifford index $k-2$, so are exceptional. 

From \cite[Thm. 3.6 and Thm. 3.7]{elms}, the Clifford dimension of any smooth 
$C \in |L|$ is $\frac{1}{2}(k+1)$ and only $\O_C(D)$ computes the Clifford index.
\end{proof}

We now recall the conjecture in \cite{elms}:

\vspace{0,3cm}

\noindent {\bf Conjecture (Eisenbud, Lange, Martens, Schreyer).}
{\it  Let $C$ be a smooth curve of Clifford dimension $r \geq 3$. Then:} 
  \begin{itemize}
  \item[(a)] {\it  $C$ has genus $g=4r-2$ and Clifford index $c=2r-3$ (and thus degree $d=g-1$),
  \item[(b)] $C$ has a unique line bundle $A$ computing $c$,
  \item[(c)] $A ^2 \iso \omega_C$ and $A$ embeds $C$ as an
    arithmetically Cohen-Macaulay curve in $\PP^r$,
  \item[(d)] $C$ is $2r$-gonal, and there is a one-dimensional family
    of pencils of degree $2r$, all of the form $|A-B|$, where $B$ is
    the divisor of $2r-3$ points of $C$.}
  \end{itemize}

\vspace{0,3cm}

In \cite{elms} the conjecture is proved for
$r \leq 9$, and in general it is proved that if $C$ satisfies (a), then it also satisfies (b)-(d).
We therefore see that the curves in the ``generalized ELMS examples'' satisfy the conjecture.

To prove Theorem 2, we use the well-known theorem of Green and Lazarsfeld.

Let $C \in |L|$ be a smooth exceptional curve on a $K3$ surface, of genus $g$, 
Clifford index $c$ and gonality $c+3$, different from a smooth plane sextic in the Donagi-Morrison example. 
Then, by Theorem 1, all smooth curves in $|L|$ have the same gonality $d=c+3$. 
As $2d-2-g=\rho(d,g,1) \leq 1$, we have $c < \lfloor \frac{g-1}{2} \rfloor$.
By \cite{gl} all the curves in $|L|$ have Clifford index $c$ and there is a line 
bundle $N$ on $S$ such that
$c=\Cliff \O_C(N)$ and (see e.g. \cite{mak3,kn,JK}) we also have that $|N|$ is base point free
$h^0(N) \geq 2$, $h^0(L-N) \geq 2$, $h^1(N)=h^1(L-N)=0$ and $N.(L-N)=c+2$.  

 By Lemma \ref{lemma:2.2} we must have
(possibly after interchanging $N$ and $L-N$)  that $L \sim 2N+\Gamma$ for a smooth rational curve $\Gamma$
satisfying $\Gamma.N=1$. In particular $c=N^2-1$ and $N^2 >0$. Therefore, the general element $D \in |N|$ is a smooth curve.

To show that we are in the ``generalized ELMS examples'' 
we have left to show that there is no line bundle $B$ on $S$ satisfying 
$0 \leq B^2 \leq N^2-1$ and $0< B.L-B^2 \leq N^2$. 

Assume such a $B$ exists. Then the numerical conditions imply $(L-B)^2 \geq N^2+3>0$
and $(L-B).L \geq 2N^2+3 >0$, so that $h^0(L-B) \geq 2$ by Riemann-Roch. Similarly $h^0(B) \geq 2$ and one therefore easily sees that $\O_C(B)$ contributes to the Clifford index of $C$, for any smooth $C \in |L|$.
Hence
\[ \Cliff C \leq \Cliff \O_C(B) \leq B.L -2(h^0(B)-1) \leq B.L-B^2-2 \leq D^2-2=c-1, \]
a contradiction.

Thus, Theorem 2 is proved. 

\begin{rem}
{\rm Note that the curves in the ``generalized ELMS examples'' have $\dim W^1_d(C)=1$ and 
$\rho(d,g,1)=0$, where $d=\gon C$ (cf. Theorem \ref{thm:zero}).}
\end{rem}

\end{document}